\documentclass[preprint,12pt]{elsarticle}

\usepackage{amssymb}
\usepackage{amsmath}
\usepackage{amssymb}
\usepackage{amsthm}
\usepackage{arydshln}
\usepackage{titletoc}
\usepackage{titlesec}
\usepackage{subfigure}
\usepackage{setspace}

\begin{document}

\begin{frontmatter}



\newtheorem{proposition}{Proposition}[section]
\newtheorem{theorem}{Theorem}[section]
\newtheorem{definition}{Definition}[section]
\newtheorem{lemma}[theorem]{Lemma}
\newtheorem{remark}{Remark}[section]
\newtheorem{example}{Example}[section]
\newtheorem{corollary}{Corollary}[section]

\title{Quaternion Gradient and Hessian\tnoteref{t1}}
\tnotetext[t1]{This work was supported by the National Natural Science Foundation of
China (No.61301202), by the Research Fund for the Doctoral Program of Higher Education of China
(No.20122304120028), and by the Fundamental Research Funds for the Central Universities.}
\author[rvt,focal]{Dongpo Xu\corref{cor1}}
\ead{dongpoxu@gmail.com}
\author[focal]{Danilo P. Mandic}
\ead{d.mandic@imperial.ac.uk}
\cortext[cor1]{Corresponding author. Tel: +44 2075946271, Fax: +44 2075946302}
\address[rvt]{College of Science, Harbin Engineering University, Harbin 150001, China}
\address[focal]{Department of Electrical and Electronic Engineering, Imperial College London, SW7 2AZ London, UK.}

\begin{abstract}
The optimization of real scalar functions of quaternion variables, such as the mean square error or array output power, underpins many practical applications.
Solutions often require the calculation of the gradient and Hessian,
however, real functions of quaternion variables are essentially non-analytic. To address this issue, we propose new definitions of quaternion gradient and Hessian, based on the novel generalized HR (GHR) calculus,
thus making possible efficient derivation of optimization algorithms directly in the quaternion field, rather than transforming the problem to the real domain, as is current practice.
In addition, unlike the existing quaternion gradients, the GHR calculus allows for the product and chain rule, and for a one-to-one correspondence of the proposed quaternion gradient and Hessian with their real counterparts.
Properties of the quaternion gradient and Hessian relevant to numerical applications are elaborated, and the results illuminate the usefulness of the GHR calculus in greatly simplifying the derivation of the quaternion least mean squares,
and in quaternion least square and Newton algorithm. The proposed gradient and Hessian are also shown to enable the same generic forms as the corresponding real- and complex-valued algorithms, further illustrating the advantages in algorithm design and evaluation.
\end{abstract}

\end{frontmatter}


\section{Introduction.}
In recent years, quaternion algebra has been successfully applied to
numerous problems in physics \cite{Girard}, computer graphics
\cite{Hanson}, signal processing and communications
\cite{Ujang11,Took09,Cyrus13,Moreno12,Bihan,Bihan14,Buchholz,Miron}. In these
applications, quaternions have allowed for a reduction in the number of
parameters and operations involved. Despite the obvious advantages, the
main obstacle for their more widespread use is that the real-valued
cost (objective) functions of quaternion variables, that naturally
arise in optimization procedures, are not analytic according to
quaternion analysis \cite{Sudbery,Deavours,Leo}. To bypass the issue
of non-existent derivatives of real functions of quaternion
variables, current optimization procedures typically rewrite a cost
function $f$ in terms of the four real components $q_a,q_b,q_c$ and
$q_d$ of the quaternion variable $q$, and take the real derivatives
with respect to these independent real variables, the so called
pseudo-derivatives \cite{Took09}. In this way, $f$ is treated as a real analytic
mapping between $\mathbb{R}^4$ and $\mathbb{R}$, however, since the
original framework is quaternion-valued, it is often awkward and
tedious to reformulate the problem in the real domain and to
calculate gradients.

To this end, the recent HR calculus \cite{Mandic11} takes the formal derivatives of $f$ with respect to
the quaternion variables and involutions $q,q^i,q^j,q^k$, and is shown to save on computational burden and
to greatly simplify gradient expressions. The HR calculus can be considered as a generalization of the complex CR calculus \cite{Brandwood,Kreutz,Wirtinger}
to the quaternion field, as the basis for the HR calculus is the use of involutions
(generalized conjugate) \cite{Ell}. However, the traditional product rule does not apply within the HR calculus because of
the non-commutativity of quaternion algebra. We here solve this issue using the novel GHR calculus,
and equip quaternion analysis with both the product rule and chain rule \cite{DPXU},
leading to the integral form of the mean value theorem and Taylor's theorem.
This also enables the derivation of optimization techniques and learning algorithms based on square error to be carried out directly
in the quaternion field, rather than transforming the problem to the real domain, and for concise and elegant algorithm development,
unlike the existing tedious pseudo-derivatives.

In this paper, a quaternion gradient and Hessian are defined based
on the GHR calculus, and in the same spirit as the complex gradient
and Hessian \cite{Bos,Kreutz}, which were instrumental for the
developments in complex domain optimization and signal processing
\cite{Sorber}. The basic relationships between the quaternion
gradient and Hessian and their real counterparts are established by
invertible linear transforms. Such mappings are shown to be very
convenient for the derivation of algorithms based on first- and
second-order Taylor series expansion in the quaternion field.
The quadrivariate
vs. quaternion relations involve redundancy, since they operate in the augmented
quaternion space $\mathbb{H}^{4N\times 1}$, and we further propose an efficient way to obtain the algorithms
that operate directly in $\mathbb{H}^{N\times 1}$. The
paper concludes with some enabling techniques for applications, such as a concise derivation of the quaternion least mean square, quaternion least square and Newton algorithms using the proposed GHR
framework.

\subsection{Quaternion Algebra}
We shall first briefly review the basic concepts in quaternion
algebra. For advanced reading on
quaternions, we refer to \cite{Ward}, as well as to \cite{FZhang97}
for several important results on matrices of quaternions.
 Quaternions are an associative but not commutative
algebra over $\mathbb{R}$, defined as
\begin{equation}\label{eq:realreprent}
\mathbb{H}=\textrm{span}\{1,i,j,k\}\triangleq \{q_a+iq_b+jq_c+kq_d\;|\;q_a,q_b,q_c,q_d\in\mathbb{R}\}
\end{equation}
where $\{1,i,j,k\}$ is a basis of $\mathbb{H}$, and the imaginary
units $i,j$ and $k$ satisfy $i^2=j^2=k^2=ijk=-1$. For a quaternion
$q=q_a+iq_b+jq_c+kq_d=S_q+V_q$, the real and vector parts are
denoted by $\mathfrak{R}(q)=q_a=S_q$ and
$\mathfrak{I}(q)=V_q=iq_b+jq_c+kq_d$. Note that the quaternion
product is noncommutative, i.e., in general for $p,q\in\mathbb{H}$, $pq\neq qp$. The conjugate of a quaternion $q$ is
$q^*=q_a-iq_b-jq_c-kq_d$, while the conjugate of the product
satisfies $(pq)^*=q^*p^*$. The modulus of a quaternion is defined as
$|q|=\sqrt{qq^*}$, and it is easy to check that $|pq|=|p||q|$. The
inverse of a quaternion $q\neq 0$ is $q^{-1}=q^*/|q|^2$, and an
important property of the inverse is $(pq)^{-1}=q^{-1}p^{-1}$. If
$|q| = 1$, we call $q$ a \textit{unit} quaternion. A quaternion $q$
is said to be \textit{pure} if $\mathfrak{R}(q)=0$, then $q^*=-q$
and $q^2=-|q|^2$. Thus, a \textit{pure unit} quaternion is a square
root of -1, such as the imaginary units $i,j$ and $k$.

Quaternions can also be written in the \textit{polar} form $
q=|q|(\cos \theta + I_q \sin \theta )$, where $I_q=Vq/|Vq|$ is a
pure unit quaternion and
$\theta=\arccos(\mathfrak{R}(q)/|q|)\in\mathbb{R}$ is the angle (or
argument) of the quaternion. We shall next introduce the quaternion
rotation and involution operations.
\begin{definition}[Quaternion Rotation \cite{Ward}]\label{def:qrot}
For any quaternion $q$, consider the transformation
\begin{equation}\label{eq:qrotat}
 q^{\mu}\triangleq \mu q \mu^{-1}
\end{equation}
where $\mu=|\mu|(\cos \theta + \hat{\mu} \sin \theta )$ is any non-zero quaternion. This transformation geometrically describes a 3-dimensional rotation of the vector part of $q$ by an angle $2\theta$ about the vector part of $\mu$.
\end{definition}

In particular, if $\mu$ in \eqref{eq:qrotat} is a pure unit
quaternion, then the quaternion rotation \eqref{eq:qrotat} becomes
the quaternion involution given in \cite{Ell}. Some important
properties of the quaternion rotation (see \cite{DPXU,Buchholz})
are:
\begin{equation}\label{pr:pqmu}
(pq)^{\mu}=p^{\mu}q^{\mu},\quad   pq=q^pp=q p^{(q^*)},\quad
\forall p, q \in\mathbb{H}
\end{equation}
\begin{equation}\label{pr:def1qmunu}
q^{\mu \nu}=(q^{\nu} )^{\mu},\quad q^{\mu*} \triangleq
(q^*)^{\mu}=(q^{\mu})^*\triangleq q^{*\mu},\quad  \forall \nu, \mu
\in\mathbb{H}
\end{equation}
These properties are the basis for the GHR calculus. Note that the
real representation in \eqref{eq:realreprent} can be easily
generalized to a general orthogonal system
$\{1,i^{\mu},j^{\mu},k^{\mu} \}$ given in \cite{DPXU,Ward}, where
the following properties hold
\begin{equation}\label{pr:prodijkmurl}
i^{\mu}i^{\mu}=j^{\mu}j^{\mu}=k^{\mu}k^{\mu}=i^{\mu}j^{\mu}k^{\mu}=-1
\end{equation}

\section{Quaternion Derivatives}
The quaternion pseudo-derivatives define or write a pseudo-derivative here are too restrictive, making the derivation of gradient based optimization algorithms of quaternion variables cumbersome and tedious \cite{Took09}.
A recent equivalent and more elegant approach is to use the HR calculus, which comprises of two
groups of derivatives: the HR-derivatives \cite{Mandic11}
\begin{align}\label{def:hrder}
{ \left(\begin{array}{cc}
\frac{\partial f}{\partial q}\\
\frac{\partial f}{\partial q^i}\\
\frac{\partial f}{\partial q^j}\\
\frac{\partial f}{\partial q^k}
\end{array}\right)}=\frac{1}{4}
\left(\begin{array}{cccc}
1 & -i  & -j & -k \\
1 & -i & j  & k\\
1 & i  & -j & k\\
1 & i  & j  & -k
\end{array}\right)
{ \left(\begin{array}{cc}
\frac{\partial f}{\partial q_a}\\
\frac{\partial f}{\partial q_b}\\
\frac{\partial f}{\partial q_c}\\
\frac{\partial f}{\partial q_d}
\end{array}\right)}
\end{align}
and the conjugate HR-derivatives

\begin{align}\label{def:hrconjder}
{ \left(\begin{array}{cc}
\frac{\partial f}{\partial q^*}\\
\frac{\partial f}{\partial q^{i*}}\\
\frac{\partial f}{\partial q^{j*}}\\
\frac{\partial f}{\partial q^{k*}}
\end{array}\right)}=\frac{1}{4}
\left(\begin{array}{cccc}
1 & i  & j & k \\
1 & i & -j  & -k\\
1 & -i  & j & -k\\
1 & -i  & -j  & k
\end{array}\right)
{ \left(\begin{array}{cc}
\frac{\partial f}{\partial q_a}\\
\frac{\partial f}{\partial q_b}\\
\frac{\partial f}{\partial q_c}\\
\frac{\partial f}{\partial q_d}
\end{array}\right)}
\end{align}
However, the traditional product rule are not valid for the HR calculus.
For example, $f(q)=|q|^2$, then $\frac{\partial |q|^2}{\partial q}=\frac{1}{2}q^*$ from \eqref{def:hrder},
but $ \frac{\partial |q|^2}{\partial q}\neq q\frac{\partial q^*}{\partial q}+\frac{\partial q}{\partial q}q^*=-\frac{1}{2}q+q^*$.
This difficulty has been solved within the framework of the GHR calculus, which equips quaternion analysis with both the novel product rule and chain rule, see \cite{DPXU} for more details.
\begin{definition}[The GHR Derivatives \cite{DPXU}]\label{def:leftghr}
Let $f:\mathbb{H}\rightarrow\mathbb{H}$. Then the left GHR derivatives of $f(q)$ with respect to $q^{\mu}$ and $q^{\mu*}$ are defined as
\begin{equation}\label{eq:leftghr}
\begin{split}
\frac{\partial f}{\partial q^{\mu}}=\frac{1}{4}\left(\frac{\partial f}{\partial q_a}-\frac{\partial f}{\partial q_b}i^{\mu}-\frac{\partial f}{\partial q_c}j^{\mu}-\frac{\partial f}{\partial q_d}k^{\mu}\right)\in\mathbb{H}\\
\frac{\partial f}{\partial q^{\mu*}}=\frac{1}{4}\left(\frac{\partial f}{\partial q_a}+\frac{\partial f}{\partial q_b}i^{\mu}+\frac{\partial f}{\partial q_c}j^{\mu}+\frac{\partial f}{\partial q_d}k^{\mu}\right)\in\mathbb{H}
\end{split}
\end{equation}
where $\mu\neq 0, \mu \in \mathbb{H}$, $\frac{\partial f}{\partial q_a},\frac{\partial f}{\partial q_b},\frac{\partial f}{\partial q_c},\frac{\partial f}{\partial q_d}\in\mathbb{H}$ are the partial derivatives of $f$ with respect to $q_a$, $q_b$, $q_c$ and $q_d$, respectively, and the set $\{1,i^{\mu},j^{\mu},k^{\mu}\}$ is a general orthogonal basis of $\mathbb{H}$.
Similarly, the right GHR derivatives are defined as
\begin{equation}
\begin{split}
\frac{\partial_r f}{\partial q^{\mu}}=\frac{1}{4}\left(\frac{\partial f}{\partial q_a}-i^{\mu}\frac{\partial f}{\partial q_b}-j^{\mu}\frac{\partial f}{\partial q_c}-k^{\mu}\frac{\partial f}{\partial q_d}\right)\in\mathbb{H}\\
\frac{\partial_r f}{\partial q^{\mu*}}=\frac{1}{4}\left(\frac{\partial f}{\partial q_a}+i^{\mu}\frac{\partial f}{\partial q_b}+j^{\mu}\frac{\partial f}{\partial q_c}+k^{\mu}\frac{\partial f}{\partial q_d}\right)\in\mathbb{H}
\end{split}
\end{equation}
\end{definition}
Some important rules of the left GHR derivatives (see \cite{DPXU})
are:
\begin{align}
&\textrm{Constant rule}:\frac{\partial (\nu f)}{\partial
q^{\mu}}=\nu\frac{\partial f}{\partial q^{\mu}}, \quad \frac{\partial (f\nu)}{\partial q^{\mu}}=\frac{\partial f}{\partial q^{\nu\mu}}\nu\label{rl:constrl1}\\
&\textrm{Product rule}: \frac{\partial (fg)}{\partial q^{\mu}}=f\frac{\partial g}{\partial q^{\mu}}+\frac{\partial (fg)}{\partial q^{g\mu}}g\label{rl:product1}\\
&\textrm{Product rule}:\frac{\partial (fg)}{\partial
q^{\mu*}}=f\frac{\partial g}{\partial q^{\mu*}}+\frac{\partial
(fg)}{\partial q^{g\mu*}}g \label{rl:product2}
\end{align}
\begin{align}
 &\textrm{Chain rule
}:\frac{\partial f(g(q))}{\partial
q^{\mu}}=\sum_{\nu\in\{1,i,j,k\}}\frac{\partial f}{\partial
g^{\nu}}\frac{\partial g^{\nu}}{\partial q^{\mu}}\label{rl:chain1}\\\
&\textrm{Chain rule }:\frac{\partial f(g(q))}{\partial
q^{\mu*}}=\sum_{\nu\in\{1,i,j,k\}}\frac{\partial f}{\partial
g^{\nu}}\frac{\partial g^{\nu}}{\partial
q^{\mu*}}\label{rl:chain2}\\
&\textrm{Rotation rule} :   \left(\frac{\partial f}{\partial
q^{\mu}}\right)^{\nu} =\frac{\partial f^{\nu}}{\partial
q^{\nu\mu}},\quad\left(\frac{\partial f}{\partial
q^{\mu*}}\right)^{\nu}
=\frac{\partial f^{\nu}}{\partial q^{\nu\mu*}}\\
&\quad \textrm{if $f$ is real } \left(\frac{\partial f}{\partial
q^{\mu}}\right)^{\nu} =\frac{\partial f}{\partial
q^{\nu\mu}},\quad\left(\frac{\partial f}{\partial
q^{\mu*}}\right)^{\nu} =\frac{\partial f}{\partial
q^{\nu\mu*}}\label{rl:realrota}\\
 &\textrm{Conjugate
rule}:\left(\frac{\partial f}{\partial q^{\mu}}\right)^*
=\frac{\partial_r f^*}{\partial q^{\mu*}},\quad\left(\frac{\partial f}{\partial q^{\mu*}}\right)^*=\frac{\partial_r f^*}{\partial q^{\mu}}\\
&\quad \textrm{if $f$ is real } \left(\frac{\partial f}{\partial
q^{\mu}}\right)^*=\frac{\partial f}{\partial q^{\mu*}},\quad
\left(\frac{\partial f}{\partial q^{\mu*}}\right)^*=\frac{\partial
f}{\partial q^{\mu}}\label{rl:realconj}
\end{align}


\begin{remark}
Observe that for $\mu\in\{1,i,j,k\}$, the HR derivatives
\eqref{def:hrder} and \eqref{def:hrconjder} are a special case of
the right GHR derivative, the latter being more concise and intuitive.  Furthermore, the GHR derivatives incorporate a novel
product rules \eqref{rl:product1}, \eqref{rl:product2} and chain
rules \eqref{rl:chain1}, \eqref{rl:chain2}, which are very
convenient for calculating the GHR derivatives.
\end{remark}

\begin{remark}
Due to the non-commutativity of quaternion products, the left GHR derivative is different from the right GHR derivative. However,
they will be equal if the function $f$ is real-valued \cite{DPXU}. In the sequel, we mainly focus on the left GHR derivative, because it has a lot of convenient properties \eqref{rl:constrl1}-\eqref{rl:realconj} that are consistent with our common sense.
\end{remark}

\section{Quaternion Gradient}

The existing notions of quaternion gradient are confusing, which has prevented systematic development of quaternion gradient based optimization.
The quaternion pseudo-gradient (also known as component-wise gradients) has been used in \cite{Took09}, however,
the calculation of the pseudo-gradient is cumbersome and tedious, making the derivation of optimization algorithms of quaternion variables very prone to errors.
An equivalent and more elegant approach is
to use the HR calculus \cite{Mandic11}, which defines the HR-gradient with respective to quaternion vector variable and its conjugate, however, the traditional product rule is not applicable to the HR-gradient.
In addition, a definition of quaternion gradient based on involutions
has been proposed in \cite{unifTook} and is termed the I-gradient,
however, the direction of the I-gradient has not be proven to follow the maximum rate of change in \cite{unifTook}, and the traditional product rule is also not valid for the I-gradient.
To this end, a quaternion gradient is here defined based on the GHR calculus, which comprises the novel product rules \eqref{rl:product1}, \eqref{rl:product2} and chain
rules \eqref{rl:chain1}, \eqref{rl:chain2}.
We consider a general case of a function $f({\bf
q}):\mathbb{H}^{N\times 1}\rightarrow\mathbb{H}$, where ${\bf
q}=(q_1,q_2,\cdots,q_N)\in \mathbb{H}^{N\times 1}$.

\begin{definition}[Quaternion Gradient]\label{def:qgrad}
The quaternion gradient and conjugate gradient of a function $f: \mathbb{H}^{N\times 1}\rightarrow \mathbb{H}$ are defined as
\begin{equation*}
\begin{split}
\nabla_{{\bf q}} f \triangleq \left(\frac{\partial f}{\partial {\bf q}}\right)^T=\left(\frac{\partial f}{\partial q_1},\cdots,\frac{\partial f}{\partial q_N}\right)^T\in\mathbb{H}^{N\times 1}\\
\nabla_{{\bf q}^*} f \triangleq \left(\frac{\partial f}{\partial
{\bf q}^*}\right)^T=\left(\frac{\partial f}{\partial
q_1^*},\cdots,\frac{\partial f}{\partial
q_N^*}\right)^T\in\mathbb{H}^{N\times 1}
\end{split}
\end{equation*}
\end{definition}
\begin{definition}[Quaternion Jacobian Matrix]\label{def:leftjaco}
If ${\bf f}: \mathbb{H}^{N\times 1}\rightarrow \mathbb{H}^{M\times1}$, then the quaternion Jacobian matrix and conjugate Jacobian matrix of ${\bf f}$ are defined as
\begin{equation*}
\begin{split}
\frac{\partial {\bf f}}{\partial {\bf q}}=\left(\begin{array}{ccc}
\frac{\partial f_1}{\partial q_1}& \cdots &\frac{\partial f_1}{\partial q_N} \\
\vdots  & \ddots & \vdots \\
\frac{\partial f_M}{\partial q_1}& \cdots &\frac{\partial
f_M}{\partial q_N}
\end{array}\right),\quad
\frac{\partial {\bf f}}{\partial {\bf q}^*}=\left(\begin{array}{ccc}
\frac{\partial f_1}{\partial q_1^*}& \cdots &\frac{\partial f_1}{\partial q_N^*} \\
\vdots  & \ddots & \vdots \\
\frac{\partial f_M}{\partial q_1^*}& \cdots &\frac{\partial
f_M}{\partial q_N^*}
\end{array}\right)\in\mathbb{H}^{M\times N}
\end{split}
\end{equation*}
\end{definition}
Note the convention that $\frac{\partial {\bf f}}{\partial {\bf q}}$ for two vectors ${\bf f}\in\mathbb{H}^{M\times 1}$ and ${\bf q}^{N\times 1}$  is a matrix
whose $(m, n)$th element is $(\partial f_m/\partial q_n)$. Thus, the dimension of $\frac{\partial {\bf f}}{\partial {\bf q}}$ is $M\times N$.
This convention will be used throughout this paper. We consider a quaternion vector ${\bf q}={\bf q}_a+i{\bf q}_b+j{\bf
q}_c+k{\bf q}_d\in\mathbb{H}^{N\times 1}$, expressed by its real
coordinate vectors ${\bf q}_a,{\bf q}_b,{\bf q}_c$ and ${\bf
q}_d\in\mathbb{R}^{N\times 1}$. Following an approach similar to that
in \cite{Mandic11,DPXU} for the case of scalar quaternions, we can now
define an augmented quaternion vector ${\bf h}\in
\mathbb{H}^{4N\times 1}$ based on its involutions, as ${\bf h}=({\bf
q}^T,{\bf q}^{iT},{\bf q}^{jT},{\bf
q}^{kT})^T\in\mathbb{H}^{4N\times 1}$. Its relationship with the
dual quadrivariate real vector in $\mathbb{R}^{4N\times 1}$ is given by \cite{Took11,Jvia10}

\begin{align}\label{eq:corrdtranmat}
{\bf h}\triangleq \left(\begin{array}{c}
{\bf q}  \\
{\bf q}^i \\
{\bf q}^j  \\
{\bf q}^k  \\
\end{array}\right)= {\bf J}{\bf r}=\left(\begin{array}{cccc}
{\bf I}_N & i{\bf I}_N  & j{\bf I}_N & k{\bf I}_N \\
{\bf I}_N & i{\bf I}_N  & -j{\bf I}_N & -k{\bf I}_N \\
{\bf I}_N & -i{\bf I}_N  & j{\bf I}_N & -k{\bf I}_N \\
{\bf I}_N & -i{\bf I}_N  & -j{\bf I}_N & k{\bf I}_N \\
\end{array}\right) \left(\begin{array}{c}
{\bf q}_a  \\
{\bf q}_b \\
{\bf q}_c  \\
{\bf q}_d
\end{array}\right)
\end{align}
where ${\bf r}=({\bf q}_a^T,{\bf q}_b^{T},{\bf q}_c^{T},{\bf
q}_d^{T})^T\in\mathbb{R}^{4N\times 1}$, ${\bf I}_N$ is the $N\times N$ identity matrix, and ${\bf J}$ is the
$4N\times 4N$ matrix in \eqref{eq:corrdtranmat}. Multiplying both
sides of \eqref{eq:corrdtranmat} by ${\bf J}^H$ and noting ${\bf
J}{\bf J}^H=4{\bf I}_{4N}$, we have
\begin{equation}\label{eq:realijksum}
\quad{\bf r}=\frac{1}{4}{\bf J}^H{\bf h}\;\in \mathbb{H}^{4N\times 1}
\end{equation}
From \eqref{eq:corrdtranmat}, a real scalar function $f({\bf
q}):\mathbb{H}^{N\times 1}\rightarrow\mathbb{R}$ can be viewed in
three equivalent forms
\begin{equation}\label{eq:threuivfrm}
f({\bf q})\quad \Leftrightarrow \quad
f({\bf q}_a,{\bf q}_b,{\bf q}_c,{\bf q}_d)\triangleq f({\bf r})\quad \Leftrightarrow \quad  f({\bf h}) \triangleq f({\bf q},{\bf q}^i,{\bf q}^j,{\bf q}^k)
\end{equation}
Since \eqref{eq:corrdtranmat} is a linear transformation and ${\bf r}$ is a real vector, it follows that
\begin{equation}\label{eq:multihr}
\frac{\partial f}{\partial {\bf h}}=\frac{\partial f}{\partial {\bf
r}}\frac{\partial {\bf r}}{\partial {\bf h}}
=\frac{1}{4}\frac{\partial f}{\partial {\bf r}}{\bf J}^H \quad
\Leftrightarrow \quad \frac{\partial f}{\partial {\bf
r}}=\frac{\partial f}{\partial {\bf h}}{\bf J},\quad f\in\mathbb{R}
\end{equation}
where $\frac{\partial f}{\partial {\bf h}}\in\mathbb{H}^{1\times 4N}$ and $\frac{\partial
f}{\partial {\bf r}}\in\mathbb{R}^{1\times 4N}$. Since $f$ and
${\bf r}$ are real-valued, we have
\begin{equation}\label{eq:realqtgradrel}
\begin{split}
\nabla_{\bf r} f &\triangleq \left(\frac{\partial f}{\partial {\bf r}}\right)^T=\left(\frac{\partial f}{\partial {\bf r}}\right)^H \\
& =\left(\frac{\partial f}{\partial {\bf h}}{\bf J}\right)^H\quad (\textrm{from}\;\eqref{eq:multihr})\\
& ={\bf J}^H\left(\frac{\partial f}{\partial {\bf h}}\right)^H\\
& ={\bf J}^H\left(\frac{\partial f}{\partial {\bf h}^*}\right)^T \quad (\textrm{from}\;\eqref{rl:realconj})\\
& ={\bf J}^H\nabla_{{\bf h}^*}f
\end{split}
\end{equation}
This shows that the real gradient $\nabla_{\bf r} f \in
\mathbb{R}^{4N\times 1}$ and the augmented quaternion gradient
$\nabla_{{\bf h}^*} f \in \mathbb{H}^{4N\times 1}$ are related by a
simple invertible linear transformation ${\bf J}^H$. From
\eqref{eq:corrdtranmat}, \eqref{eq:threuivfrm} and
\eqref{eq:realqtgradrel}, we can now state that each of the following
equations represents a necessary and sufficient condition for the
existence of stationary points of a real-valued function $f$
\begin{equation}
\frac{\partial f}{\partial {\bf q}}={\bf
0}\;\Leftrightarrow\;\frac{\partial f}{\partial {\bf q}^*}={\bf 0}\;
\Leftrightarrow \; \frac{\partial f}{\partial {\bf r}}={\bf 0} \;
\Leftrightarrow \; \frac{\partial f}{\partial {\bf h}}={\bf 0}\;
\Leftrightarrow \; \frac{\partial f}{\partial {\bf h}^*}={\bf 0}
\end{equation}

\subsection{Quaternion Gradient Descent Algorithm}\label{sec:qgda1}
Gradient descent (also known as steepest descent) is a first-order optimization algorithm, which finds a local minimum of a function by taking steps proportional to the negative of the gradient of the function at the current point.
For a real scalar function $f({\bf q}):\mathbb{H}^{N\times 1}\rightarrow\mathbb{R}$, it can also be viewed as $f({\bf r}):\mathbb{R}^{4N\times 1}\rightarrow\mathbb{R}$ from \eqref{eq:threuivfrm}, so that the quadrivariate real gradient descent update rule can be given by \cite{Haykin,Nocedal}
\begin{equation}\label{eq:realgdrule}
\Delta{\bf r}=-\alpha \nabla_{\bf r} f,\quad {\bf r}\in
\mathbb{R}^{4N\times 1}
\end{equation}
where $\Delta{\bf r}$ denotes a small increment in ${\bf r}$ and
$\alpha\in\mathbb{R}^+$ is the step size. Using
\eqref{eq:corrdtranmat}, \eqref{eq:realqtgradrel}, and
\eqref{eq:realgdrule}, we now obtain
\begin{equation}\label{eq:qtgdrule}
\Delta {\bf h}={\bf J}\Delta {\bf r}=-\alpha {\bf J} \nabla_{\bf r}
f =-\alpha {\bf J} {\bf J}^H\nabla_{{\bf h}^*}
f=-4\alpha\nabla_{{\bf h}^*} f
\end{equation}
From \eqref{eq:corrdtranmat}, we have ${\bf h}=({\bf q}^T,{\bf
q}^{iT},{\bf q}^{jT},{\bf q}^{kT})^T$, so that \eqref{eq:qtgdrule} can be
rewritten as
\begin{equation}\label{eq:deltaugquatgrad}
\Delta {\bf h}=\left(\begin{array}{c}
\Delta{\bf q}  \\
\Delta{\bf q}^i \\
\Delta{\bf q}^j  \\
\Delta{\bf q}^k  \\
\end{array}\right)
=-4\alpha\left(\begin{array}{c}
\nabla_{{\bf q}^*} f  \\
\nabla_{{\bf q}^{i*}} f \\
\nabla_{{\bf q}^{j*}} f  \\
\nabla_{{\bf q}^{k*}} f  \\
\end{array}\right)
\end{equation}
This gives the quaternion gradient descent (QGD) update rule in the form
\begin{equation}\label{eq:qgd1}
\Delta{\bf q}=-4\alpha \nabla_{{\bf q}^*}
f=-4\alpha\left(\frac{\partial f}{\partial {\bf q^*}}\right)^T=-4\alpha\left(\frac{\partial f}{\partial {\bf q}}\right)^H,\quad
f\in\mathbb{R}
\end{equation}

\begin{remark}
From \eqref{eq:qgd1}, the quaternion gradient of a real-valued
scalar function $f$ with respect to a quaternion vector ${\bf q}$ is
equal to $\nabla_{{\bf q}^*} f =\left(\frac{\partial f}{\partial
{\bf q}^*}\right)^T=\left(\frac{\partial f}{\partial {\bf
q}}\right)^H$, but not $\nabla_{{\bf q}} f$. This result is a
generalization of the complex CR calculus given in
\cite{Brandwood,Kreutz}, and makes possible compact derivation of
learning algorithms in $\mathbb{H}$.
\end{remark}

\section{Quaternion Hessian}
Since a formal derivative of a function $f:\mathbb{H}\rightarrow\mathbb{H}$ is (wherever it exists) again a function from $\mathbb{H}$ to $\mathbb{H}$, it makes sense to take the GHR derivative of a GHR derivative, that is, a higher order GHR derivative. We shall
consider second order quaternion derivatives of the form
\begin{equation}
\begin{split}
\frac{\partial^2 f}{\partial q^{\mu} \partial q^{\nu}}=\frac{\partial }{\partial q^{\mu}}\left(\frac{\partial f}{\partial q^{\nu}}\right), \quad
\frac{\partial^2 f}{\partial q^{\mu*} \partial q^{\nu*}}=\frac{\partial }{\partial q^{\mu*} }\left(\frac{\partial f}{\partial q^{\nu*} }\right)\\
\frac{\partial^2 f}{\partial q^{\mu} \partial q^{\nu*}}=\frac{\partial }{\partial q^{\mu}}\left(\frac{\partial f}{\partial q^{\nu*}}\right), \quad
\frac{\partial^2 f}{\partial q^{\mu*} \partial q^{\nu}}=\frac{\partial }{\partial q^{\mu*} }\left(\frac{\partial f}{\partial q^{\nu} }\right)
\end{split}\quad \forall\mu,\nu \in \mathbb{H}
\end{equation}
The second order cross-derivatives are in general not identical \cite{DPXU}, that is
\begin{equation}
\frac{\partial^2 f}{\partial q^{\mu} \partial q^{\nu}}\neq\frac{\partial^2 f}{\partial q^{\nu} \partial q^{\mu}}
\end{equation}
However, the second order GHR derivatives have a commutative property \cite{Sudbery}
\begin{equation}\label{pr:realmumuconj}
16\frac{\partial^2 f}{\partial q^{\mu} \partial
q^{\mu*}}=16\frac{\partial^2 f}{\partial q^{\mu*} \partial q^{\mu}}
= \frac{\partial^2 f}{\partial q^2_a}+\frac{\partial^2 f}{\partial
q^2_b} +\frac{\partial^2 f}{\partial q^2_c}+\frac{\partial^2
f}{\partial q^2_d}
\end{equation}
If $f$ is a real-valued function, then the conjugate rule of the second order GHR derivatives is given by \cite{DPXU}
\begin{equation}\label{pr:conjsecderirel}
\begin{split}
\left(\frac{\partial^2 f}{\partial q^{\mu} \partial q^{\nu}}\right)^*=\frac{\partial^2 f}{\partial q^{\nu*} \partial q^{\mu*}},\quad
\left(\frac{\partial^2 f}{\partial q^{\mu*} \partial q^{\nu*}}\right)^*=\frac{\partial^2 f}{\partial q^{\nu} \partial q^{\mu}}\\
\left(\frac{\partial^2 f}{\partial q^{\mu} \partial q^{\nu*}}\right)^*=\frac{\partial^2 f}{\partial q^{\nu} \partial q^{\mu*}},\quad
\left(\frac{\partial^2 f}{\partial q^{\mu*} \partial q^{\nu}}\right)^*=\frac{\partial^2 f}{\partial q^{\nu*} \partial q^{\mu}}\\
\end{split}
\end{equation}

\begin{definition}[Quaternion Hessian Matrix]\label{def:lefthessian}
Let $f: \mathbb{H}^{N\times 1}\rightarrow \mathbb{H}$, then the quaternion Hessian matrix of the mapping $f$ is defined as
\begin{equation*}
{\bf H}_{{\bf q}{\bf q}} \triangleq \frac{\partial }{\partial {\bf q}}\left(\frac{\partial f}{\partial {\bf q}}\right)^T=\left(\begin{array}{ccc}
\frac{\partial^2 f}{\partial q_1\partial q_1}& \cdots &\frac{\partial^2 f}{\partial q_N\partial q_1} \\
\vdots  & \ddots & \vdots \\
\frac{\partial^2 f}{\partial q_1\partial q_N}& \cdots &\frac{\partial^2 f}{\partial q_N\partial q_N} \\
\end{array}\right)\in \mathbb{H}^{N\times N}
\end{equation*}

\begin{equation*}
{\bf H}_{{\bf q}{\bf q}^*} \triangleq \frac{\partial }{\partial {\bf q}}\left(\frac{\partial f}{\partial {\bf q}^*}\right)^T=\left(\begin{array}{ccc}
\frac{\partial^2 f}{\partial q_1\partial q_1^*}& \cdots &\frac{\partial^2 f}{\partial q_N\partial q_1^*} \\
\vdots  & \ddots & \vdots \\
\frac{\partial^2 f}{\partial q_1\partial q_N^*}& \cdots &\frac{\partial^2 f}{\partial q_N\partial q_N^*} \\
\end{array}\right)\in\mathbb{H}^{N\times N}
\end{equation*}
\end{definition}

Using \eqref{pr:realmumuconj} and \eqref{pr:conjsecderirel}, it then follows that ${\bf H}_{{\bf q}{\bf q}^*}$ is Hermitian,
so that ${\bf H}^{H}_{{\bf q}{\bf q}^*} ={\bf H}_{{\bf q}{\bf q}^*}$. Then, up to second order,
the Taylor series expansion (TSE) of the real scalar function $f({\bf q}):\mathbb{H}^{N\times 1}\rightarrow \mathbb{R}$ viewed as an analytic function $f({\bf r}):\mathbb{R}^{4N\times 1}\rightarrow \mathbb{R}$ of the vector ${\bf r}\in \mathbb{R}^{4N\times 1}$ from \eqref{eq:threuivfrm}, is given by \cite{Haykin,Nocedal}
\begin{equation}\label{eq:tsereal}
f({\bf r}+\Delta {\bf r})=f({\bf r})+\frac{\partial f}{\partial {\bf r}}\Delta {\bf r}+\frac{1}{2}\Delta {\bf r}^T{\bf H_{rr}}\Delta {\bf r}+\textrm{h.o.t}.
\end{equation}
where ${\bf H}_{{\bf rr}} \triangleq \frac{\partial }{\partial {\bf r}}\left(\frac{\partial f}{\partial {\bf r}}\right)^T\in\mathbb{R}^{4N\times 4N}$ is the real symmetric Hessian matrix, ${\bf H}^T_{{\bf rr}}={\bf H}_{{\bf rr}}$, and h.o.t. are the higher order terms. From \eqref{eq:corrdtranmat} and \eqref{eq:multihr}, the first order term in the augmented quaternion space is calculated as
\begin{equation}\label{eq:tsefirst}
\frac{\partial f}{\partial {\bf r}}\Delta {\bf r}=\frac{1}{4}\frac{\partial f}{\partial {\bf h}}{\bf J}{\bf J}^H\Delta{\bf h}=\frac{\partial f}{\partial {\bf h}}\Delta{\bf h}
\end{equation}
Noting from \eqref{eq:corrdtranmat} that ${\bf h}=({\bf q}^T,{\bf q}^{iT},{\bf q}^{jT},{\bf q}^{kT})^T$, we can now expand the first order term in \eqref{eq:tsereal} as follows
\begin{equation}\label{eq:realtsefirst}
\begin{split}
&\frac{\partial f}{\partial {\bf r}}\Delta {\bf r}=\frac{\partial f}{\partial {\bf h}}\Delta{\bf h}\quad (\textrm{from}\;\eqref{eq:tsefirst})\\
&=\frac{\partial f}{\partial {\bf q}}\Delta{\bf q} +\frac{\partial
f}{\partial {\bf q}^i}\Delta{\bf q}^i +\frac{\partial f}{\partial
{\bf q}^j}\Delta{\bf q}^j
+\frac{\partial f}{\partial {\bf q}^k}\Delta{\bf q}^k\quad (\textrm{from \eqref{eq:corrdtranmat}})\\
&=\frac{\partial f}{\partial {\bf q}}\Delta{\bf q}
+\left(\frac{\partial f}{\partial {\bf q}}\Delta{\bf q}\right)^i
+\left(\frac{\partial f}{\partial {\bf q}}\Delta{\bf q}\right)^j
+\left(\frac{\partial f}{\partial {\bf q}}\Delta{\bf q}\right)^k\quad (\textrm{from \eqref{rl:realrota}})\\
&=4\mathfrak{R}\left\{\frac{\partial f}{\partial {\bf q}}\Delta{\bf q}\right\}
\quad (\textrm{from \eqref{eq:realijksum}})
\end{split}
\end{equation}
Now, we shall consider the $4N\times 4N$ augmented quaternion Hessian matrix
\begin{equation}\label{eq:augquathess}
{\bf H}_{{\bf h}{\bf h}^*} \triangleq \frac{\partial }{\partial {\bf h}}\left(\frac{\partial f}{\partial {\bf h}^*}\right)^T
=\left(\begin{array}{cccc}
{\bf H}_{{\bf q}{\bf q}^*} & {\bf H}_{{\bf q}^i{\bf q}^*}  & {\bf H}_{{\bf q}^j{\bf q}^*} & {\bf H}_{{\bf q}^k{\bf q}^*} \\
{\bf H}_{{\bf q}{\bf q}^{i*}} & {\bf H}_{{\bf q}^i{\bf q}^{i*}}  & {\bf H}_{{\bf q}^j{\bf q}^{i*}} & {\bf H}_{{\bf q}^k{\bf q}^{i*}} \\
{\bf H}_{{\bf q}{\bf q}^{j*}} & {\bf H}_{{\bf q}^i{\bf q}^{j*}}  & {\bf H}_{{\bf q}^j{\bf q}^{j*}} & {\bf H}_{{\bf q}^k{\bf q}^{j*}} \\
{\bf H}_{{\bf q}{\bf q}^{k*}} & {\bf H}_{{\bf q}^i{\bf q}^{k*}}  & {\bf H}_{{\bf q}^j{\bf q}^{k*}} & {\bf H}_{{\bf q}^k{\bf q}^{k*}} \\
\end{array}\right)
\end{equation}
Its equivalence with ${\bf H}_{{\bf rr}}\in\mathbb{R}^{4N\times 4N}$ can be established as
\begin{equation}\label{eq:realqtnhessrel}
\begin{split}
{\bf H}_{{\bf rr}} &= \frac{\partial }{\partial {\bf r}}\left(\frac{\partial f}{\partial {\bf r}}\right)^T\\
&=\frac{\partial }{\partial {\bf r}}\left(\frac{\partial f}{\partial {\bf r}}\right)^H \quad (\textrm{since $f$ is  real-valued})\\
& =\frac{\partial }{\partial {\bf r}}\left(\frac{\partial f}{\partial {\bf h}}{\bf J}\right)^H\quad (\textrm{from}\;\eqref{eq:multihr})\\
& =\frac{\partial }{\partial {\bf r}}\left\{{\bf J}^H\left(\frac{\partial f}{\partial {\bf h}}\right)^H\right\}\\
& =\frac{\partial }{\partial {\bf h}}\left\{{\bf J}^H\left(\frac{\partial f}{\partial {\bf h}}\right)^H\right\}{\bf J} \quad (\textrm{from}\;\eqref{eq:multihr})\\
& ={\bf J}^H\frac{\partial }{\partial {\bf h}}\left(\frac{\partial f}{\partial {\bf h}^*}\right)^T{\bf J} \quad (\textrm{from}\;\eqref{rl:realconj})\\
& ={\bf J}^H{\bf H}_{{\bf h}{\bf h}^*}{\bf J}
\end{split}
\end{equation}
Note here that the Hermitian operator in \eqref{eq:realqtnhessrel} can not be replaced with the transpose operator, because quaternion matrices $({\bf AB})^T\neq {\bf B}^T{\bf A}^T$.
Recalling that the Hessian ${\bf H}_{{\bf rr}}$ is a real symmetric matrix, it is evident from \eqref{eq:realqtnhessrel} that
${\bf H}_{{\bf h}{\bf h}^*}$ is Hermitian, that is, ${\bf H}^H_{{\bf h}{\bf h}^*}={\bf H}_{{\bf h}{\bf h}^*}$. Subsequently, for the second order term of \eqref{eq:tsereal}, we have
\begin{equation}\label{eq:tsesecond}
\begin{split}
\frac{1}{2}\Delta {\bf r}^T{\bf H_{rr}}\Delta {\bf r}&=\frac{1}{2}\Delta {\bf r}^H{\bf H_{rr}}\Delta {\bf r}\quad (\textrm{since ${\bf r}$ is real-valued})\\
&=\frac{1}{2}\Delta {\bf r}^H{\bf J}^H{\bf H}_{{\bf h}{\bf h}^*}{\bf J} \Delta {\bf r}\quad (\textrm{from }\;\eqref{eq:realqtnhessrel})\\
&=\frac{1}{2}\Delta {\bf h}^H{\bf H}_{{\bf h}{\bf h}^*}\Delta {\bf h}\quad (\textrm{from }\;\eqref{eq:corrdtranmat})
\end{split}
\end{equation}
Now, from \eqref{eq:corrdtranmat} and \eqref{eq:tsesecond}, the second order term in \eqref{eq:tsereal} can be expanded as

\begin{equation}\label{eq:realtsesecond}
\begin{split}
\frac{1}{2}\Delta {\bf r}^T{\bf H_{rr}}\Delta {\bf r}&=\frac{1}{2}\Delta {\bf h}^H{\bf H}_{{\bf h}{\bf h}^*}\Delta {\bf h}\\
&=\frac{1}{2}\sum_{\mu,\nu\in\{1,i,j,k\}}\left(\Delta {\bf q}^{\nu}\right)^H{\bf H}_{{\bf q}^{\mu}{\bf q}^{\nu*}}\Delta {\bf q}^{\mu}
\quad (\textrm{from }\;\eqref{eq:augquathess})\\
&=2\sum_{\mu\in\{1,i,j,k\}}\mathfrak{R}\left(\Delta {\bf q}^H{\bf H}_{{\bf q}^{\mu}{\bf q}^*}\Delta {\bf q}^{\mu}\right)
\quad (\textrm{from \eqref{eq:realijksum} and \eqref{rl:realrota}})\\
\end{split}
\end{equation}
Thus, using \eqref{eq:threuivfrm}, \eqref{eq:tsereal}, \eqref{eq:tsefirst} and \eqref{eq:tsesecond}, the TSE expansion in $\mathbb{H}^{4N}$ (augmented
TSE) up to the second term can be expressed as
\begin{equation}\label{eq:tseaugment}
f({\bf h}+\Delta {\bf h})=f({\bf h})+\frac{\partial f}{\partial {\bf h}}\Delta {\bf h}+\frac{1}{2}\Delta {\bf h}^H{\bf H}_{{\bf h}{\bf h}^*}\Delta {\bf h}+\textrm{h.o.t.}
\end{equation}
Finally, combining the expansions given in \eqref{eq:realtsefirst} and \eqref{eq:realtsesecond} yields the TSE expressed directly in $\mathbb{H}^N$, given by
\begin{equation}\label{eq:tseregular}
\begin{split}
f({\bf q}+\Delta {\bf q})&=f({\bf q})+4\mathfrak{R}\left(\frac{\partial f}{\partial {\bf
q}}\Delta {\bf
q}\right)\\
&\quad+2\sum_{\mu\in\{1,i,j,k\}}\mathfrak{R}\left(\Delta {\bf
q}^H{\bf H}_{{\bf q}^{\mu}{\bf q}^*}\Delta {\bf
q}^{\mu}\right)+\textrm{h.o.t}.
\end{split}
\end{equation}

Using relation \eqref{eq:realqtnhessrel} and noting that $\frac{1}{4}{\bf J}{\bf J}^H={\bf I}_{4N}$, we have
\begin{equation}\label{eq:equivrrhhhess}
{\bf H}_{{\bf h}{\bf h}^*}-\lambda{\bf I}_{4N}=\frac{1}{16}{\bf J}\left({\bf H}_{\bf rr}-4\lambda{\bf I}_{4N}\right){\bf J}^H
\end{equation}

\begin{remark}
Since the real scalar function $f$ is non-analytic, the TSE in \eqref{eq:tseaugment} and \eqref{eq:tseregular} is always augmented, due to the presence of the terms $\Delta {\bf q}^{\mu}$, $\mu\in\{i,j,k\}$.
This is in contrast to the complex TSE for an analytic function.
\end{remark}

\begin{remark}
Equation \eqref{eq:equivrrhhhess} illustrates that the eigenvalues of the quadrivariate real Hessian ${\bf H}_{\bf rr}$ are quadruple of those of the augmented quaternion Hessian ${\bf H}_{{\bf h}{\bf h}^*}$. An important consequence is that the augmented quaternion Hessian ${\bf H}_{{\bf h}{\bf h}^*}$ and quadrivariate real Hessian ${\bf H}_{\bf rr}$ have the same positive definiteness properties and condition number. This result is important
in numerical applications using the Hessian such as quaternion Newton minimization.
\end{remark}

\section{Application Examples}
Quaternion gradient and Hessian are particularly useful for analytic and numerical solutions of
quaternion valued parameter estimation problems.

\subsection{Quaternion Newton Algorithm}
Newton's method is a second order optimization
method which makes use of the Hessian matrix.
This method often has a better convergence than the gradient descent method, but
can be very expensive to calculate and store the Hessian matrix.
For a real scalar criterion $f({\bf q}):\mathbb{H}^{N\times 1}\rightarrow\mathbb{R}$, which from \eqref{eq:threuivfrm} can be viewed as $f({\bf r}):\mathbb{R}^{4N\times 1}\rightarrow\mathbb{R}$
, the Newton iteration step $\Delta {\bf r}$ for the minimisation
of the function $f({\bf r})$ with respect to its real parameters ${\bf r}=({\bf q}_a^T,{\bf q}_b^{T},{\bf q}_c^{T},{\bf q}_d^{T})^T$ is described by \cite{Haykin,Nocedal}
\begin{equation}\label{eq:realnewton}
{\bf H}_{\bf rr}\Delta {\bf r}=-\nabla_{\bf r} f, \quad {\bf r}\in \mathbb{R}^{4N\times 1}
\end{equation}
where $\nabla_{\bf r} f$ is the real gradient defined by \eqref{eq:realqtgradrel}. From \eqref{eq:corrdtranmat}, \eqref{eq:realqtgradrel} \eqref{eq:realqtnhessrel} and \eqref{eq:realnewton}, it then follows that
\begin{equation}
{\bf J}^H{\bf H}_{{\bf h}{\bf h}^*}\Delta {\bf h}={\bf J}^H{\bf H}_{{\bf h}{\bf h}^*}{\bf J}\Delta {\bf r}={\bf H}_{\bf rr}\Delta {\bf r}= -\nabla_{\bf r} f=-{\bf J}^H\nabla_{{\bf h}^*} f
\end{equation}
Thus, the Newton method in the augmented quaternion domain can be formulated as
\begin{equation}\label{eq:augquatnewton}
{\bf H}_{{\bf h}{\bf h}^*}\Delta {\bf h}=-\nabla_{{\bf h}^*} f
\end{equation}
where ${\bf H}_{{\bf h}{\bf h}^*}$ is the augmented quaternion Hessian defined by \eqref{eq:augquathess}. Since ${\bf h}=({\bf q}^T,{\bf q}^{iT},{\bf q}^{jT},{\bf q}^{kT})^T$, we can rewrite \eqref{eq:augquatnewton} as
\begin{equation}
\left(\begin{array}{cccc}
{\bf H}_{{\bf q}{\bf q}^*} & {\bf H}_{{\bf q}^i{\bf q}^*}  & {\bf H}_{{\bf q}^j{\bf q}^*} & {\bf H}_{{\bf q}^k{\bf q}^*} \\
{\bf H}_{{\bf q}{\bf q}^{i*}} & {\bf H}_{{\bf q}^i{\bf q}^{i*}}  & {\bf H}_{{\bf q}^j{\bf q}^{i*}} & {\bf H}_{{\bf q}^k{\bf q}^{i*}} \\
{\bf H}_{{\bf q}{\bf q}^{j*}} & {\bf H}_{{\bf q}^i{\bf q}^{j*}}  & {\bf H}_{{\bf q}^j{\bf q}^{j*}} & {\bf H}_{{\bf q}^k{\bf q}^{j*}} \\
{\bf H}_{{\bf q}{\bf q}^{k*}} & {\bf H}_{{\bf q}^i{\bf q}^{k*}}  & {\bf H}_{{\bf q}^j{\bf q}^{k*}} & {\bf H}_{{\bf q}^k{\bf q}^{k*}} \\
\end{array}\right)
\left(\begin{array}{c}
\Delta{\bf q}  \\
\Delta{\bf q}^i \\
\Delta{\bf q}^j  \\
\Delta{\bf q}^k  \\
\end{array}\right)
=-\left(\begin{array}{c}
\nabla_{{\bf q}^*} f  \\
\nabla_{{\bf q}^{i*}} f \\
\nabla_{{\bf q}^{j*}} f  \\
\nabla_{{\bf q}^{k*}} f  \\
\end{array}\right)
\end{equation}
If ${\bf H}_{{\bf h}{\bf h}^*}$ (equivalently, ${\bf H}_{{\bf r}{\bf r}}$ in \eqref{eq:equivrrhhhess}) is positive definite,
then using the Banachiewicz inversion formula for the inverse of a nonsingular partitioned matrix \cite{FZhang}, yields
\begin{equation}\label{eq:invaugnewtonform}
\left(\begin{array}{c}
\Delta{\bf q}  \\
\Delta{\bf q}^i \\
\Delta{\bf q}^j  \\
\Delta{\bf q}^k  \\
\end{array}\right)
=-\left(\begin{array}{cc}
{\bf H}^{-1}_{{\bf q}{\bf q}^*}+{\bf L}{\bf T}^{-1}{\bf U} &-{\bf L}{\bf T}^{-1} \\
-{\bf T}^{-1}{\bf U} &{\bf T}^{-1}
\end{array}\right)\left(\begin{array}{c}
\nabla_{{\bf q}^*} f  \\
\nabla_{{\bf q}^{i*}} f \\
\nabla_{{\bf q}^{j*}} f  \\
\nabla_{{\bf q}^{k*}} f  \\
\end{array}\right)
\end{equation}
where ${\bf T}=\left({\bf H}_{{\bf h}{\bf h}^*}/{\bf H}_{{\bf q}{\bf q}^*}\right)$ is the Schur complement \cite{FZhang}
of ${\bf H}_{{\bf q}{\bf q}^*}$ in ${\bf H}_{{\bf h}{\bf h}^*}$, and
\begin{equation}\label{eq:LU}
{\bf L}={\bf H}^{-1}_{{\bf q}{\bf q}^*}\left(\begin{array}{c}
{\bf H}_{{\bf q}{\bf q}^{i*}} \\
{\bf H}_{{\bf q}{\bf q}^{j*}}  \\
{\bf H}_{{\bf q}{\bf q}^{k*}}  \\
\end{array}\right)^H,\quad
{\bf U}=\left(\begin{array}{c}
{\bf H}_{{\bf q}{\bf q}^{i*}} \\
{\bf H}_{{\bf q}{\bf q}^{j*}}  \\
{\bf H}_{{\bf q}{\bf q}^{k*}}  \\
\end{array}\right){\bf H}^{-1}_{{\bf q}{\bf q}^*}
\end{equation}
The invertibility of the Schur complement $\bf T$ follows from the positive definiteness of ${\bf H}_{{\bf h}{\bf h}^*}$, so that the quaternion Newton update rule is given by
\begin{equation}\label{eq:invnewtonform}
\Delta{\bf q} = -{\bf H}^{-1}_{{\bf q}{\bf q}^*}\nabla_{{\bf q}^*} f
+ {\bf L}{\bf T}^{-1}
\left(\begin{array}{c}
-{\bf H}_{{\bf q}{\bf q}^{i*}}{\bf H}^{-1}_{{\bf q}{\bf q}^*}\nabla_{{\bf q}^{*}} f+\nabla_{{\bf q}^{i*}} f \\
-{\bf H}_{{\bf q}{\bf q}^{j*}}{\bf H}^{-1}_{{\bf q}{\bf q}^*}\nabla_{{\bf q}^{*}} f+\nabla_{{\bf q}^{j*}} f  \\
-{\bf H}_{{\bf q}{\bf q}^{k*}}{\bf H}^{-1}_{{\bf q}{\bf q}^*}\nabla_{{\bf q}^{*}} f+\nabla_{{\bf q}^{k*}} f  \\
\end{array}\right)
\end{equation}
A substantial simplification can be introduced by avoiding the computation of  the inverse of the Schur complement ${\bf T}$, so that
the quaternion Newton (QN) method in \eqref{eq:invnewtonform} can be approximated as
\begin{equation}\label{eq:appnewton}
\Delta{\bf q}\thickapprox -{\bf H}^{-1}_{{\bf q}{\bf q}^*}\nabla_{{\bf q}^*} f
\end{equation}
\begin{remark}
Note that the redundancy in \eqref{eq:invaugnewtonform} has been removed and the resulting QN algorithm \eqref{eq:invnewtonform} and approximated QN algorithm \eqref{eq:appnewton} operate directly in $\mathbb{H}^{N\times 1}$.
For the interested reader, we leave the estimation problem of the upper bound of the approximation error between \eqref{eq:invnewtonform} and \eqref{eq:appnewton}.
\end{remark}


\subsection{Quaternion Least Mean Square}
In this subsection, we derive the quaternion least mean square (QLMS) algorithm in \cite{Took09,Mandic11} using the GHR calculus.
For convenience, the QLMS derivation applied component-wise can be found in Appendix.
Within QLMS, the cost function to be minimized is a real-valued function
\begin{equation}
J(n)=|e(n)|^2=e^*(n)e(n)
\end{equation}
where
\begin{equation}\label{eq:qlmy}
e(n)=d(n)-{\bf w}^T(n){\bf x}(n),\quad e^*(n)=d^*(n)-{\bf x}^H(n){\bf w}^*(n)
\end{equation}
$d(n)\in \mathbb{H}$ and ${\bf w}(n),{\bf x}(n)\in \mathbb{H}^{N\times 1}$. From \eqref{eq:qgd1}, the weight update of QLMS is then given by
\begin{equation}\label{eq:qlmsllearnrule}
{\bf w}(n+1)-{\bf w}(n)=-\alpha\nabla_{{\bf w}^*} J(n)=-\alpha\left(\frac{\partial J(n)}{\partial {\bf w}^*}\right)^T=-\alpha\left(\frac{\partial J(n)}{\partial {\bf w}}\right)^H
\end{equation}
where $\alpha$ is the step size and the negative gradient $-\nabla_{{\bf w}^*} J(n)$ defines the direction of gradient descent in \eqref{eq:qgd1}.
By using the product rule in \eqref{rl:product1}, the gradient is calculated as
\begin{equation}\label{eq:qlmsdJdwconj}
\frac{\partial J(n)}{\partial {\bf w}}=e^*(n)\frac{\partial e(n)}{\partial {\bf w}}+\frac{\partial e^*(n)}{\partial {\bf w}^{e(n)}}e(n)
\end{equation}
The above two derivatives now become
\begin{equation}\label{eq:qlmsdedwconj}
\begin{split}
\frac{\partial e(n)}{\partial {\bf w}}&=-\frac{\partial \left({\bf w}^T(n){\bf x}(n)\right)}{\partial {\bf w}}=-\mathfrak{R}({\bf x}(n))\\
\frac{\partial e^*(n)}{\partial {\bf w}^{e(n)}}e(n)&=-\frac{\partial \left({\bf x}^H(n){\bf w}^*(n)\right)}{\partial {\bf w}^{e(n)}}e(n)
=\frac{1}{2}{\bf x}^H(n)e^*(n)
\end{split}
\end{equation}
where the terms $\frac{\partial (q\nu)}{\partial q}$ and $\frac{\partial (\omega q^*)}{\partial q^{\mu}}{\mu}$ are given in \cite{DPXU}, and are used in the last equalities in the expressions above. Substituting \eqref{eq:qlmsdedwconj} into \eqref{eq:qlmsdJdwconj} yields
\begin{equation}\label{eq:qlmsdJdwconj1b}
\begin{split}
\frac{\partial J(n)}{\partial {\bf w}}&=-e^*(n)\mathfrak{R}({\bf x}(n))+\frac{1}{2}{\bf x}^H(n)e^*(n)\\
&=\left(\frac{1}{2}{\bf x}^H(n)-\mathfrak{R}({\bf x}(n))\right)e^*(n)=-\frac{1}{2}{\bf x}^T(n)e^*(n)
\end{split}
\end{equation}
Finally, the update of the adaptive weight vector of QLMS becomes
\begin{equation}\label{eq:qlmsllearnrule1b}
{\bf w}(n+1)={\bf w}(n)+\frac{1}{2}\alpha\, e(n){\bf x}^*(n)
\end{equation}
where the constant $\frac{1}{2}$ can be absorbed into $\alpha$.

\begin{remark}
From \eqref{eq:qlmsllearnrule1b} and \eqref{eq:compwiseLMS}, we can see that the QLMS derived using the GHR calculus is exactly the same as that using the pseudo-gradient,
however, the derivation of component-wise gradient in Appendix is too cumbersome and tedious.
The equality of \eqref{eq:qlmsllearnrule1b} and \eqref{eq:compwiseLMS} also provides a theoretical support for the gradient descent method in \eqref{eq:qgd1}.
Note that if we start from $e(n)=d(n)-{\bf w}^H(n){\bf x}(n)$, the final update rule of QLMS would become ${\bf w}(n+1)={\bf w}(n)+\alpha\, {\bf x}(n)e^*(n)$.
The QLMS algorithm in \eqref{eq:qlmsllearnrule1b} is a therefore generalization of complex
LMS \cite{Widrow} to the case of quaternion vector.
\end{remark}

\begin{remark}
The QLMS algorithm \eqref{eq:qlmsllearnrule1b} is different from the original QLMS \cite{Took09} based on
componentwise gradients, the HR-QLMS \cite{Mandic11} based on the HR-gradient, and the I-QLMS \cite{unifTook} based on the I-gradient.
The difference with the original QLMS arises due to the rigorous use of the non-commutativity of quaternion product in \eqref{eq:wabcdJn} and \eqref{eq:compwiseLMS}.
The difference with the HR-QLMS and I-QLMS is due to the rigorous use of the novel product rule in \eqref{eq:qlmsdJdwconj}.
\end{remark}

\subsection{Quaternion Least Squares}
The quaternion least squares (QLS) problem can be
formulated as: Given ${\bf A}\in \mathbb{H}^{M\times N}$ and ${\bf
b}\in \mathbb{H}^{M\times 1}$ with $M\geq N$, find ${\bf q}\in \mathbb{H}^{N\times 1}$ such that
 the error or residual of the overdetermined linear system of equations
\begin{equation}\label{eq:qlsproblem}
F({\bf q})=\|{\bf b}-{\bf A}{\bf q}\|^2=\left({\bf b}-{\bf A}{\bf q}\right)^H\left({\bf b}-{\bf A}{\bf q}\right)
\end{equation}
is minimized. Upon taking the GHR derivative of $F({\bf q})$ using the product rule \eqref{rl:product1} and setting to zero, we have
\begin{equation}\label{eq:lqsgrad}
\begin{split}
\frac{\partial F({\bf q})}{\partial {\bf q}}&=\frac{\partial \left({\bf b}^H{\bf b}-{\bf b}^H{\bf Aq}-{\bf q}^H{\bf A}^H{\bf b}+{\bf q}^H{\bf A}^H{\bf Aq} \right)}{\partial {\bf q}}\\
&=-{\bf b}^H{\bf A}+\frac{1}{2}\left({\bf A}^H{\bf b}\right)^H+{\bf q}^H{\bf A}^H{\bf A}-\frac{1}{2}\left({\bf A}^H{\bf Aq}\right)^H \\
&=-\frac{1}{2}\left({\bf b-Aq}\right)^H{\bf A}={\bf 0}
\end{split}
\end{equation}
Thus, the overdetermined system in
\eqref{eq:qlsproblem} reduces to an ($N\times N$) linear system
(normal equation)
\begin{equation}
{\bf A}^H{\bf Aq}={\bf A}^H{\bf b}
\end{equation}
In this way, for a nonsingular ${\bf A}^H{\bf A}$, the unique solution of \eqref{eq:qlsproblem} is given by
\begin{equation}\label{equation}
{\bf q}=\left({\bf A}^H{\bf A}\right)^{-1}{\bf A}^H{\bf b}={\bf A}^+{\bf b}
\end{equation}
where ${\bf A}^+$ denotes the Moore-Penrose generalized
inverse in \cite{Yuan2013}. From \eqref{eq:qgd1} and \eqref{eq:lqsgrad}, the gradient of $F$ can be expressed as
\begin{equation}
\nabla_{{\bf q}^*} F =\left(\frac{\partial F({\bf q})}{\partial {\bf
q}}\right)^H=-\frac{1}{2}{\bf A}^H\left({\bf b-Aq}\right)
\end{equation}
while using the constant rule \eqref{rl:constrl1}, the quaternion Hessian matrices of $F$ is calculated as
\begin{equation}\label{eq:qlshess}
{\bf H}_{{\bf q}{\bf q}^*} =\frac{\partial }{\partial {\bf
q}}\left(\frac{\partial F({\bf q})}{\partial {\bf
q}}\right)^H=\frac{1}{2}{\bf A}^H{\bf A},\quad {\bf H}_{{\bf q}^i{\bf q}^*}={\bf H}_{{\bf q}^j{\bf q}^*}={\bf H}_{{\bf q}^k{\bf q}^*}={\bf 0}
\end{equation}
In this case, the matrix ${\bf L}$ in \eqref{eq:LU} and \eqref{eq:invnewtonform} becomes zero, so the approximation error between \eqref{eq:invnewtonform} and \eqref{eq:appnewton} vanishes.
\begin{remark}
The optimal quaternion solution \eqref{equation} is formally equivalent to the real representation matrix (17) in \cite{Yuan2013},
however, the real and complex representation methods of quaternion matrix equations are too cumbersome and difficult to use.
Such quaternion matrix equations can be solved directly in the quaternion field using the GHR calculus.
\end{remark}

\section{Conclusions}
A new formulation for the quaternion gradient
and Hessian of a smooth real function of
quaternion variables has been proposed based on the GHR calculus. It has been shown that the so obtained quaternion
gradient and Hessian and their real counterparts are related by
simple linear transforms. The quaternion gradient descent and quaternion Newton algorithm have
been derived and shown to perform all computations
directly in the quaternion field, without the need to increase
the problem dimensionality. The GHR calculus thus resolves the long standing problems of quaternion analyticity, product and chain rule,
and greatly simplifies the derivation of first- and second-order iterative optimisation procedures.
The proposed framework has been shown to serve as a basis for generic extensions of real- and complex-valued optimization solutions.
Apart from the gradient descent, least squares and Newton algorithms addressed in this work, the conjugate-gradient and quasi-Newton methods
are also readily obtained.

\section*{Acknowledgments.}
We thank Dr. Cyrus Jahanchahi and Dr. Clive Cheong
Took for fruitful discussions on the HR calculus.

\section*{Appendix: The QLMS derivation applied component-wise}
From Definitions \ref{def:leftghr} and \ref{def:qgrad},
the weight update of QLMS can be written componentwise as
\begin{equation}\label{eq:realgrad}
\begin{split}
&{\bf w}(n+1)-{\bf w}(n)=-\alpha\nabla_{{\bf w}^*} J(n)\\
&=-\frac{1}{4}\alpha\left(\nabla_{{\bf w}_a} J(n)+\nabla_{{\bf w}_b} J(n)i+\nabla_{{\bf w}_c} J(n)j+\nabla_{{\bf w}_d} J(n)k\right)
\end{split}
\end{equation}
where $\alpha$ is the step size and the negative gradient $-\nabla_{{\bf w}^*} J(n)$ defines the direction of gradient descent in \eqref{eq:qgd1}.
By using the traditional product rule, the subgradients in \eqref{eq:realgrad} is therefore calculated by
\begin{equation}\label{eq:subgradsss}
\begin{split}
\nabla_{{\bf w}_a} J(n)=e^*(n)(\nabla_{{\bf w}_a} e(n))+(\nabla_{{\bf w}_a} e^*(n))e(n)\\
\nabla_{{\bf w}_b} J(n)=e^*(n)(\nabla_{{\bf w}_b} e(n))+(\nabla_{{\bf w}_b} e^*(n))e(n)\\
\nabla_{{\bf w}_a} J(n)=e^*(n)(\nabla_{{\bf w}_c} e(n))+(\nabla_{{\bf w}_c} e^*(n))e(n)\\
\nabla_{{\bf w}_d} J(n)=e^*(n)(\nabla_{{\bf w}_d} e(n))+(\nabla_{{\bf w}_d} e^*(n))e(n)
\end{split}
\end{equation}
where the validity of traditional product rule is owing to the real valued nature of ${\bf w}_a,{\bf w}_b,{\bf w}_c$ and ${\bf w}_d$.
We can now  calculate the following subgradient
\begin{align}\label{eq:subgrad}
\nabla_{{\bf w}_a} e(n)&=-\nabla_{{\bf w}_a} ({\bf w}^T(n){\bf x}(n))=-\nabla_{{\bf w}_a} \left(({\bf w}^T_a+{\bf w}^T_bi+{\bf w}^T_cj+{\bf w}^T_dk){\bf x}(n)\right)\nonumber\\
&=-\nabla_{{\bf w}_a} \left(({\bf w}^T_a){\bf x}(n)\right)=-{\bf x}(n)
\end{align}
and similarly
\begin{equation}
\nabla_{{\bf w}_b} e(n)=-i{\bf x}(n),\quad\nabla_{{\bf w}_c} e(n)=-j{\bf x}(n),\quad\nabla_{{\bf w}_d} e(n)=-k{\bf x}(n)
\end{equation}
Following on \eqref{eq:subgrad}, the subgradients of $e^*(n)$ in \eqref{eq:subgradsss} can be expressed as
\begin{align}\label{eq:subgradconj}
\nabla_{{\bf w}_a} e^*(n)&=-\nabla_{{\bf w}_a} ({\bf x}^H(n){\bf w}^*(n))=-\nabla_{{\bf w}_a} \left({\bf x}^H(n)({\bf w}_a-{\bf w}_bi-{\bf w}_cj-{\bf w}_dk)\right)\nonumber\\
&=-\nabla_{{\bf w}_a} \left({\bf x}^H(n){\bf w}_a\right)=-{\bf x}^*(n)
\end{align}
In a similar manner, we have
\begin{equation}\label{eq:subgradconj2}
\nabla_{{\bf w}_b} e^*(n)={\bf x}^*(n)i,\quad\nabla_{{\bf w}_c} e^*(n)={\bf x}^*(n)j,\quad\nabla_{{\bf w}_d} e^*(n)={\bf x}^*(n)k
\end{equation}
Now, upon substituting \eqref{eq:subgrad}-\eqref{eq:subgradconj2} to \eqref{eq:subgradsss}, we have
\begin{equation}\label{eq:wabcdJn}
\begin{split}
\nabla_{{\bf w}_a} J(n)=-e^*(n){\bf x}(n)-{\bf x}^*(n)e(n)\\
\nabla_{{\bf w}_b} J(n)=-e^*(n)i{\bf x}(n)+{\bf x}^*(n)ie(n)\\
\nabla_{{\bf w}_c} J(n)=-e^*(n)j{\bf x}(n)+{\bf x}^*(n)je(n)\\
\nabla_{{\bf w}_d} J(n)=-e^*(n)k{\bf x}(n)+{\bf x}^*(n)ke(n)
\end{split}
\end{equation}
Finally, substituting \eqref{eq:wabcdJn} to \eqref{eq:realgrad}, we arrive at the expression of the QLMS in the form
\begin{equation}\label{eq:compwiseLMS}
\begin{split}
{\bf w}(n+1)-{\bf w}(n)&=-\alpha\nabla_{{\bf w}^*} J(n)\\
&=\frac{1}{4}\alpha e^*(n)\left({\bf x}(n)+i\,{\bf x}(n)\,i+j\,{\bf x}(n)\,j+k\,{\bf x}(n)\,k\right)\\
&\quad+\frac{1}{4}\alpha {\bf x}^*(n)\left(e(n)-i\,e(n)\,i-j\,e(n)\,j-k\,e(n)\,k\right)\\
&=-\frac{1}{2}\alpha e^*(n){\bf x}^*(n)+\alpha {\bf x}^*(n)\mathfrak{R}(e(n))\\
&=\alpha\left(-\frac{1}{2} e^*(n)+\mathfrak{R}(e(n))\right) {\bf x}^*(n)\\
&=\frac{1}{2}\alpha e(n) {\bf x}^*(n)
\end{split}
\end{equation}

\label{}





\bibliographystyle{model1a-num-names}
\bibliography{<your-bib-database>}



\end{document}